\documentclass{amsart}
\usepackage{amssymb}
\usepackage{mathrsfs}
\usepackage{graphicx}
\usepackage{comment}

\def\S{\mathcal S}

\newcommand{\nni}{\noindent}
\newcommand{\be}{\begin{equation}}
\newcommand{\ee}{\end{equation}}
\newcommand{\ba}{\begin{align}}
\newcommand{\ea}{\end{align}}

\newcommand{\abs}[1]{\lvert#1\rvert}

\newtheorem{example}{Example}[section]
\newtheorem{theorem}{Theorem}[section]

{\begin{list}{}{%
\settowidth{\labelwidth}{\textsf{{\it #1.}}}%
\setlength{\labelsep}{4mm}%
\setlength{\leftmargin}{\labelwidth}%
\addtolength{\leftmargin}{\labelsep}%
}}%
{\end{list}}

\def\beq{\begin{equation}}\def\enq{\end{equation}}

{\begin{list}{}{%
\settowidth{\labelwidth}{\textsf{{\it #1.}}}%
\setlength{\labelsep}{2mm}%
\setlength{\leftmargin}{\labelwidth}%
\addtolength{\leftmargin}{\labelsep}%
\addtolength{\leftmargin}{4mm}%
\setlength{\itemsep}{6pt}%
\setlength{\listparindent}{0pt}%
\setlength{\topsep}{3pt}%
}}%

\title[Integer group determinants ]{The integer group determinants for the non-abelian groups of order 18}

\author[B. Paudel]{Bishnu Paudel}
\address{ Department of Mathematics\\
         Kansas State University\\
         Manhattan, KS 66506, USA}
\email{%humbertb@ksu.edu,
 bpaudel@ksu.edu, pinner@math.ksu.edu}

\author[C. Pinner]{Chris Pinner}

\date{\today}
\begin{document}

\begin{abstract}
We obtain a complete description of the integer group determinants for the non-abelian groups of order 18.
%SmallGroup(18,4), 
%the generalized dihedral group for $\mathbb Z_3 \times \mathbb Z_3$. 

\end{abstract}

\maketitle

\section{Introduction} 
A problem of  Olga Taussky-Todd  \cite{TausskyTodd} is to determine the values taken by the group determinant when the entries are integers.
Here we shall think of the integer group
determinant as being defined on the elements in the group ring $\mathbb Z[G]$: 
$$ \mathcal{D}_G\left( \sum_{g\in G} a_g g \right):=\det\left( a_{gh^{-1}}\right),$$
where $g$ indexes rows  and $h$ columns. We write 
$ \mathcal{S}(G)$
for the set of integer group determinants for $G$. We use $\mathbb Z_n$ and $D_n$ 
to denote the cyclic and dihedral groups of order $n$.

 A complete description of the integer group determinants was  obtained for $\mathbb Z_p$ and $\mathbb Z_{2p}$ in  \cite{Newman1,Laquer}, for $D_{2p}$ and $D_{4p}$ in \cite{dihedral},
and for the remaining  groups of order at most 15 in  \cite{smallgps,bishnu1}.
Recently this has been extended to the 14 groups of order 16 (the 5 abelian groups in \cite{Yamaguchi1,Yamaguchi2,Yamaguchi3,Yamaguchi4,Yamaguchi5}), and  the  9 non-abelian groups in \cite{dihedral,ZnxH,Q16,Yamaguchi6,Yamaguchi7, Humb, Yamaguchi8, Humb2}).

There are  five groups of order 18. The integer group determinants for the dihedral group
were given in \cite{dihedral}:
$$ \mathcal{S}(D_{18})=\{  2^a 3^b m\; \mid \; a=0 \text{ or } a\geq 2, \; b=0 \text{ or } b\geq 5, \; \gcd(m,6)=1   \}. $$
This is SmallGroup(18,1), using  the group identification from
GAP's small group library.

In this paper we determine  the integer group determinants for the remaining two non-abelian groups. The  first
is the generalized dihedral group for $\mathbb Z_3 \times \mathbb Z_3$
$$ G_{18,4}=   \langle X,Y,Z\; | \; X^3=Y^3=Z^2=1,\; XY=YX, \; XZ=ZX^{-1},\; YZ=ZY^{-1}\rangle. $$
This is SmallGroup(18,4).

\begin{theorem} \label{G1} The integer group determinants for SmallGroup(18,4) coprime to 6
are the $m\equiv \pm 1$ mod 18. 

The determinants divisible by 6 are the $2^2\cdot 3^9 m$, $m\in\mathbb Z$. 

The determinants divisible by 2 but not 3 are the $2^2(9m\pm 2)$, $m\in \mathbb Z.$

The  determinants divisible by 3 but not 2 are the $3^9(2m+1)$, $m\in \mathbb Z.$

\end{theorem}
We give the proof of Theorem \ref{G1}  in Section 2.

The other group is 
$$\mathbb Z_3 \times  D_6 =\langle X, Y, Z\; | \; X^3=Y^3=Z^2=1, \; XZ=ZX^{-1},\; YX=XY, \; YZ=ZY\rangle. $$
This is SmallGroup(18,3).

\begin{theorem} \label{G2} The integer group determinants for $\mathbb Z_3 \times D_6$ coprime to 6
are the $m\equiv \pm 1$ mod 18. 

The determinants divisible by 6 are the $2^2\cdot 3^6 m$, $m\in\mathbb Z$. 

The determinants divisible by 2 but not 3 are the $2^2(9m\pm 2)$, $m\in \mathbb Z.$

The  determinants divisible by 3 but not 2 are the $3^6(2m+1)$, $m\in \mathbb Z.$

\end{theorem}
We give the proof of Theorem \ref{G2} in Section 3.

\subsection{The integer group determinants of a subgroup}

It remains an open problem whether $\mathcal{S}(G)\subseteq \mathcal{S}(H)$ for all subgroups $H$ of $G$. This has been proved for abelian $H$, see \cite[Theorem 1.4]{Yamaguchi1}. In \cite{ZnxH} we also showed that $\mathcal{S}(\mathbb Z_2 \times H)\subseteq \mathcal{S}(H).$ 

Notice $G_{18,4}$ and $\mathbb Z_3 \times D_6$ both have subgroups $\mathbb Z_3 \times \mathbb Z_3=\langle X,Y\rangle $,
where  
$$ \mathcal{S}(\mathbb Z_3\times \mathbb Z_3)=\{9m\pm 1 \text{ and } 3^6m \; \mid \; m\in \mathbb Z\},$$
from \cite{smallgps}. Indeed we will use the fact that we  can explicitly write our determinants as a $\mathbb Z_3\times \mathbb Z_3$ determinant. They both have  a non-abelian subgroup
  $D_6=\langle X, Z\rangle$, and  satisfy $\mathcal{S}(G)\subseteq \mathcal{S}(D_6)$, 
where  from \cite{dihedral}
$$\mathcal{S}(D_{2p})=\{ 2^ap^bm\; \mid \; \gcd(m,2p)=1,\; a=0 \text{ or } a\geq 2,\; b=0 \text{ or } b\geq 3\},$$
(with  $\langle Y, Z\rangle$ giving another $D_6$ for $G_{18,4}$  and $\mathbb Z_6$
for $\mathbb Z_3 \times D_6$).

The general affine group of order 20 has
$\mathcal{S}(GA(1,5)) \subset \mathcal{S}(D_{10}) $
%=\{ 2^a\cdot 5^d m\; \mid\; \gcd(m,10)=1,\; a=0 \text{ or } a\geq 2,\;a=0 \text{ or } a\geq 3\}, $$
(it can be shown that  a  $GA(1,5)$ determinant $2^a 5^bm$, $\gcd(m,10)=1,$  has $a=0$ or $a\geq 4$ and $b=0$ or $b\geq 5$).
Hence a counter example will require $G$ to have at least order 24.

\section{The Proof of Theorem \ref{G1}}

\subsection{A formula for the generalised dihedral group determinant}\label{formula}

Recall for an abelian group $H,$ the generalized dihedral group is the semidirect product of $H$ with the cyclic group of order two, with the non-identity element acting as the inverse map on $H$:
$$ G= \langle H, Z\; \mid \; Z^2=1, \; \;Zh=h^{-1}Z \; \;\forall h \in H \rangle. $$

Frobenius \cite{Frob} observed that the group determinant can be factored 

$$ \mathcal{D}_G\left( \sum_{g\in G} a_g g \right)=\prod_{\rho\in \hat{G}} \det \left(\sum_{g\in G} a_g \rho(g) \right)^{\deg \rho}, $$
where $\hat{G}$ denotes a complete set of non-isomorphic irreducible group  representations  for the group (see for example \cite{Conrad} or \cite{book}).

We suppose that 
$H=\mathbb Z_{m_1} \times \cdots \times \mathbb Z_{m_r}$ and write
$$H=\langle X_1,\ldots ,X_r\; \mid \; X_1^{m_1}=\ldots =X_r^{m_r}=1,\;\; X_iX_j=X_jX_i\;\;  \forall i,j\rangle,$$ 
so that
an element of $\mathbb Z[G]$ can be written
$$F= f(X_1,\ldots ,X_r)+ Zg(X_1,\ldots ,X_r) $$
where the $f(x_1,\ldots ,x_r),g(x_1,\ldots , x_r)$ are integer polys of degree at most $m_i-1$ in the $x_i$. The representations on $H$ are the $|H|$ characters where the $\chi(x_i)$
are the various $m_i$th roots of unity. Notice that
$$ \rho(Z)=\begin{pmatrix} 0 & 1 \\ 1 & 0 \end{pmatrix},\;\; \quad  \rho(h)=\begin{pmatrix} \chi(h) & 0 \\ 0 & \chi(h)^{-1}\end{pmatrix}, $$
satisfies the group relations and gives a representation on $G$ with
\begin{align*}  \det \left(\sum_{g\in G} a_g \rho(g)\right) &  =\det \begin{pmatrix} f(\chi(X_1),\ldots ,\chi(X_r)) & g\left(\chi(X_1^{-1}),\ldots ,\chi(X_r^{-1})\right) \\ g(\chi(X_1),\ldots ,\chi(X_r))  & f\left(\chi(X_1^{-1}),\ldots ,\chi(X_r^{-1}) \right)\end{pmatrix}\\
& = f(\chi(X_1),\ldots ,\chi(X_r))  f(\chi(X_1^{-1}),\ldots ,\chi(X_r^{-1})) \\
& \quad - g(\chi(X_1),\ldots ,\chi(X_r)) g(\chi(X_1^{-1}),\ldots ,\chi(X_r^{-1})). \end{align*}
Notice that these will correspond to the product of two characters if all the $\chi(X_i)=1$ or 
$-1$ (if $m_i$ is even) and degree two representations otherwise (where the same factor
occurs for $\chi$ and $\chi^{-1}$). Hence
$$\mathcal{D}_G(F) = \prod_{x_1^{m_1}=1,\ldots ,x_r^{m_r}=1} f(x_1,\ldots ,x_r)f(x_1^{-1},\ldots ,x_{r}^{-1}) -  g(x_1,\ldots ,x_r)g(x_1^{-1},\ldots ,x_{r}^{-1}). $$
Note, this is the $H$ determinant for 
$$  f(X_1,\ldots ,X_r)f(X_1^{-1},\ldots ,X_{r}^{-1}) -  g(X_1,\ldots ,X_r)g(X_1^{-1},\ldots ,X_{r}^{-1})\in \mathbb Z[H].$$

\subsection{Determinant for SmallGroup(18,4)}

We have $H=\mathbb Z_3\times \mathbb Z_3$ and can write
\be \label{groupring}  F(X,Y,Z)=f(X,Y) +Zg(X,Y) \ee
where
$$f(x,y)=a_0(x)+ya_1(x)+y^2a_2(x),\quad g(x,y)=b_0(x)+yb_1(x)+y^2b_2(x), $$
with the $a_i(x),b_i(x)\in \mathbb Z[x]$ at most quadratic.

For $G=$SmallGroup(18,4) we get
\begin{align*} \mathcal{D}_G(F)  & = \prod_{x^3=1,y^3=1} f(x,y)f(x^{-1},y^{-1})- g(x,y)g(x^{-1},y^{-1}). 
\end{align*}

With $\omega:=e^{2\pi i/3}$ denoting a primitive cuberoot of unity, we can write this as
 $$ \mathcal{D}_G(F)= A B_1^2B_2^2B_3^2B_4^2, $$
with integers
\begin{align*} 
A: & = f(1,1)^2-g(1,1)^2=\left(f(1,1)+g(1,1)\right)(f(1,1)-g(1,1)),\\
B_1: & = f(1,\omega)f(1,\omega^2)- g(1,\omega)g(1,\omega^2),\\
B_2: & = f(\omega,1)f(\omega^2,1)- g(\omega,1)g(\omega^2,1),\\
B_3: & = f(\omega,\omega)f(\omega^2,\omega^2)- g(\omega,\omega)g(\omega^2,\omega^2),\\
B_4: & = f(\omega,\omega^2)f(\omega^2,\omega)- g(\omega,\omega^2)g(\omega^2,\omega).
\end{align*}

\subsection{SmallGroup(18,4) determinants must be of the stated form}

Notice that if $2\mid A$ then $4\mid A$ (since it is the product of two terms of the same parity) while the other terms $B_1,\ldots ,B_4$ are all squared. Hence if $2\mid \mathcal{D}_G(F)$ then $2^2\mid \mathcal{D}_G(F)$.

Notice that the $B_i\equiv A$ mod $1-\omega$ and so, being integers, are congruent mod 
3. Therefore, if $3\mid \mathcal{D}_G(F),$ then $3\mid A,B_1,\ldots , B_4,$ and $3^9\mid \mathcal{D}_G(F)$.

Hence the multiples of 6 must be multiples of $2^2\cdot 3^9$ and the multiples of 3 but not 2 must be the odd multiples of $3^9$. All these can be obtained.

Since they are $\mathbb Z_3\times \mathbb Z_3$ determinants we know that the values
coprime to 3 must be $\pm 1$ mod 9.
Thus the odd values must be $\pm 1$ mod 18 and the even values must be of the form
$2^2k$ with $k\equiv \pm 2$ mod 9. Again, we can achieve all these.

\subsection{Achieving the SmallGroup(18,4) determinants}
Writing
$$ h(x,y):=(x^2+x+1)(y^2+y+1), $$
we obtain the $1+18m$ from taking
$$ f=1+mh,\quad g=mh, $$
and the $-1-18m$ by switching $f$ and $g$. 

We achieve the $2^2(2+9m)$ with
$$f=1+x+x^2+y-y^2x^2+mh,\quad g=1+x+yx-y^2(x+x^2)+mh, $$
and the $2^2(-2-9m)$  from switching $f$ and $g$.

For the multiples of 6, we get $2^2\cdot 3^{10} m$ from 
$$f=1+x-yx^2-y^2(1+x)-mh,\;\;\; g=1+y-y^2(1+x+x^2)+mh, $$
and $2^2\cdot 3^9 (1+3m)$ from
$$ f= 1+x+x^2+y(1+x)-y^2x^2+mh,\quad g=(1+x^2)(1+y-y^2)+mh, $$
with $2^2\cdot 3^9 (-1-3m)$ by switching $f$ and $g$.

For the odd multiples of 3, we get $3^{10}(1+2m)$ from
$$ f= 1+x+x^2+y(1+x+x^2)+y^2(1-x-x^2)+mh,\quad g=1+x+x^2+y(1+x)+y^2(1-x-x^2)+mh, $$
and $3^9(1+6m)$ from
$$f=1+x+x^2+y(1+x)- y^2(1+x+x^2)+mh,\quad g=1+x+x^2+y(1+x-x^2)-y^2(1+x+x^2)+mh, $$
with $3^9(-1-6m)$ from switching $f$ and $g$. \qed

\section{Proof of Theorem \ref{G2}}

\subsection{The group determinant for $G=\mathbb Z_3 \times D_6$}

From the group presentation  we can again write elemements of $\mathbb Z[G]$ in the form \eqref{groupring}.
From \cite{ZnxH} we can write the $G$ determinant as the product of three $D_6=\langle X,Z\rangle$  determinants, namely of 
$$F(X,1,Z),\quad F(X,\omega,Z), \quad F(X,\omega^2,Z), \quad \omega:=e^{2\pi i/3}. $$
The coefficients of these are in $\mathbb Z[\omega]$ rather than $\mathbb Z$. From \cite{dihedral} this gives us
$$  \mathcal{D}_G(F)= \prod_{y^3=1} \prod_{x^3=1} f(x,y)f(x^{-1},y)-g(x,y)g(x^{-1},y). $$
Notice, this is the $\mathbb Z_3 \times \mathbb Z_3=\langle X,Y\rangle$ measure 
of $f(X,Y)f(X^{-1},Y)-g(X,Y)g(X^{-1},Y)$.
Pairing conjugates we can write this in the form
$$ \mathcal{D}_G(F) = A_1A_2^2   A_3 A_4^2$$
with the integers
\begin{align*}
A_1 &=f(1,1)^2-g(1,1)^2=(f(1,1)+g(1,1))(f(1,1)-g(1,1)),\\
A_2 &=f(\omega,1)f(\omega^2,1)-g(\omega,1)g(\omega^2,1),\\
A_3 &=(f(1,\omega)^2-g(1,\omega)^2 )(f(1,\omega^2)^2-g(1,\omega^2)^2)\\
 & =\abs{f(1,\omega)+g(1,\omega)}^2\abs{f(1,\omega)-g(1,\omega)}^2,\\
A_4 & = (f(\omega,\omega)f(\omega^2,\omega) -  g(\omega,\omega)g(\omega^2,\omega))(f(\omega^2,\omega^2)f(\omega,\omega^2) -  g(\omega^2,\omega^2)g(\omega,\omega^2)).
 \end{align*}

\subsection{The $\mathbb Z_3 \times D_6$ determinants must be of the stated form}

Since $A_1$ and $A_3$ are both products of two integers with the same parity
and $A_2$ and $A_4$ are squared, the even determinants must be multiples of 4.

Now $A_2\equiv A_1$ mod 3 and $A_3,A_4\equiv A_1^2$ mod 3. Hence, if the determinant is a multiple of 3, we must have $3\mid A_1,A_2,A_3,A_4$ and $3^6\mid \mathcal{D}_G(F).$ Since it is a $\mathbb Z_3 \times \mathbb Z_3 $ determinant, the determinants coprime to 3 must be $\pm 1$ mod 9.

Hence the odd determinants coprime to 3 must be $\pm 1$ mod 18 and the even determinants coprime to 3 must be of the form $4k$ with $k\equiv \pm 2$ mod 9.
The odd multiples of $3$ must be odd multiples of $3^6$ and the even multiples of 3 
are multiples of $2^2\cdot 3^6$. We show that all these can be achieved.

\subsection{Achieving the determinants}

We write $h(x,y):=(1+x+x^2)(1+y+y^2)$.

For the values coprime to 6, we get the $1+18m$ from
$$ f=1+mh, \quad g=0+mh, $$
the $-1-18m$ from switching $f$ and $g$. 

For the multiples of 6, we get  the $2^2\cdot 3^6(1+3m)$ from 
$$ f=(1+x)+y+y^2+mh,\quad g=1+y+mh, $$
the $2^2\cdot 3^6(-1-3m)$ by switching $f$ and $g$, and the $2^2\cdot 3^7m$ from
$$ f= (1+x-x^2)+y(1-x^2)+y^2(1-x^2) +mh,\quad g= (1-x^2)+y(1+x-x^2)+y^2(x-x^2) -mh.$$

For the odd multiples of 3, we get $3^6(1+6m)$ from
$$ f=1+y+mh,\quad g=1+mh, $$
the $3^6(-1-6m)$ by switching $f$ and $g$, and the $3^7(2m+1)$ from
$$ f=(1+x)+y(1+x)+y^2(1+x-x^2)  + mh, \quad g= (1+x)+y(1+x-x^2)+y^2  +mh. $$

For the even determinants coprime to 3, we get $2^2(2+9m)$ from
$$ f=(1+x-x^2)+y x^2+ y^2 x^2+mh,\quad g= (1+x-x^2)+y(-x+x^2)+y^2(-1+x^2)+mh, $$
and $2^2(-2-9m)$ by switching $f$ and $g$. \qed

\end{document}